\documentclass[12pt]{amsart}
\usepackage{amssymb}

\newtheorem{theorem}{Theorem}[section]
\newtheorem{lemma}{Lemma}[section]

\newtheorem{remark}{Remark}[section]

\newtheorem{corollary}{Corollary}[section]

\def\a{\omega_{1}}
\def\b{\omega_{2}}
\def\c{\omega_{3}}
\def\la{\lambda}

\begin{document}
\title{Dynamically convex Finsler metrics and $J$-holomorphic embedding
of asymptotic cylinders}
\author{Adam Harris}
\address{School of Mathematics, Statistics and Computer Sciences,
University of New England,
Armidale NSW 2351,
Australia}
\email{adamh@turing.une.edu.au}
\author{Gabriel P. Paternain}
 \address{ Department of Pure Mathematics and Mathematical Statistics,
University of Cambridge,
Cambridge CB3 0WB, England}
 \email {g.p.paternain@dpmms.cam.ac.uk}

\date{}

\begin{abstract} We explore the relationship between contact forms on
$\mathbb S^3$ defined by Finsler metrics on $\mathbb S^2$ and the theory
developed by H. Hofer, K. Wysocki and E. Zehnder in \cite{HWZ,HWZ1}.
We show that a Finsler metric on $\mathbb S^2$ with curvature $K\geq 1$
and with all geodesic loops of length $>\pi$ is dynamically convex and hence
it has either two or infinitely many closed geodesics.
We also explain how to explicitly construct $J$-holomorphic embeddings
of cylinders asymptotic to Reeb orbits of contact structures arising
from Finsler metrics on $\mathbb S^2$ with $K=1$ thus complementing the 
results obtained in \cite{HW}.
\end{abstract}


\maketitle

\vspace{.2in}

%\newpage 

\section{Introduction}
A contact form $\lambda$ on a closed, connected, oriented manifold of
odd dimension is said to be {\em dynamically convex} if the Conley-Zehnder
index of any contractible periodic orbit of the Reeb vector field is at least three.
This notion generalises the case of convex hypersurfaces in ${\Bbb R}
^{4}$ on which the contact form is simply the restriction of the standard
one-form $\lambda_{0}$ from the ambient space. When the hypersurface is
${\Bbb S}^{3}$, the Reeb flow induces a fibration by orbits of constant
period, namely the Hopf fibration. If the hypersurface is an ellipsoid
${\Bbb E}_{p,q}$ defined by the equation $p|z|^{2} + q|w|^{2} = 1$
($z$, $w$ denoting coordinates of ${\Bbb C}^{2}$) then an interesting
dichotomy arises between the dynamics of the Reeb flow for $\frac{p}{q}$
rational or irrational. In the rational case there are still infinitely
many periodic orbits, though not all are of the same minimal period, while 
in the
irrational case there are just two periodic orbits. This dichotomy extends
to any strictly convex hypersurface, since
it was shown by Hofer, Wysocki and Zehnder \cite{HWZ} that the Reeb vector 
field $X$ of any contact form satisfying the dynamic convexity condition on 
the three--sphere has either two or infinitely many periodic orbits. A central
role in their theory is played by the existence of finite--energy 
$J$--holomorphic embeddings
of the cylinder in the ``symplectisation'' of a contact manifold, from which
a periodic orbit is realised as the asymptotic limit of an end of the embedded 
cylinder. 

For the standard example associated with the Hopf fibration on
${\Bbb S}^{3}$ there is a direct correspondence between the asymptotic
limit of finite--energy $J$--holomorphic embeddings in the symplectisation 
on ${\Bbb S}^{3}\times{\Bbb R}$ and germs of plane algebroid curves centred
at the origin of a sufficiently small ball in ${\Bbb C}^{2}$, via their 
``link'' with the boundary--sphere. It may then be asked whether a more general
correspondence can be established between plane algebroid curve--germs and
$J$--holomorphic embeddings in the symplectisation of other contact forms
$\lambda$ on ${\Bbb S}^{3}$. It should be noted that a similar question has 
been studied extensively in the context of pseudoholomorphic curves in a 
symplectic manifold (cf., e.g., \cite{Sik}) though the question there is local 
rather
than asymptotic in nature. In \cite{HW} it was shown that near periodic 
orbits of
``elliptic type'' (cf. {\em locally recurrent} orbits, \cite{HW}) such that 
the 
partial almost complex endomorphism $j$, compatible with the contact structure,
is Reeb--invariant (i.e., the Lie derivative ${\mathcal L}_{X}j = 0$), finite--energy
$J$--holomorphic embeddings $\psi$ of the cylinder may be represented
holomorphically in a suitable tubular coordinate neighbourhood of the orbit.
After transformation of the cylinder to a punctured disc $D\setminus\{0\}$ via 
the choice of a complex coordinate $z$, the periodic component of $\psi$ 
naturally subdivides $D$ into ``quasi--sectors'', the number of these
being determined by the topological degree of $\psi$ restricted to the circle
$|z| = c$,  as $c$ approaches
zero (i.e., as the axial coordinate of the cylinder approaches infinity),
and being referred to as the ``charge'' of the mapping. In the direction
spanned by a disc $\Delta$, transversal to the orbit, $\psi$ is represented by
holomorphic functions defined on each of the quasi--sectors of $D$,
which are hinged together along common edges by the return map $\alpha$,
induced on $\Delta$ by the Reeb flow. With respect to $(\lambda, j)$ satisfying
the above criteria and for a given charge $n$, each finite--energy $J$--
holomorphic embedding of the cylinder in ${\Bbb S}^{3}\times{\Bbb R}$
therefore gives rise to a collection of $n$ holomorphic functions, whose 
continuity at adjacent boundaries of their domains is mediated by $\alpha$ 
(cf. \cite[Theorem 1]{HW}).

For any contact structure $(\lambda, j)$ such that
${\mathcal L}_{X}j = 0$ and $\alpha$ is the identity (e.g., the standard 
$\lambda_{0}$ restricted to ${\Bbb S}^{3}$) it is possible to move explicitly 
back and forth between finite--energy $J$--holomorphic embeddings of the 
cylinder in the tubular neighbourhood of a periodic orbit and algebraic
curve--germs at the origin in ${\Bbb C}^{2}$ of the form $(z^{n}, F(z))$.
A similarly explicit correspondence was obtained in \cite{HW} for contact 
structures
on ${\Bbb S}^{3}$ diffeomorphically equivalent to the restriction of 
$\lambda_{0}$
to a rational ellipsoid (on which $\alpha$ is a non--trivial rational rotation
near each of two exceptional orbits). Implicit in these
examples is the fact that the contact form has rotational symmetry along the
axis of a suitable tubular coordinate neighbourhood of the orbit, thus 
allowing the map $\psi$ conversely to be determined (up to local diffeomorphic
equivalence) by the charge $n$ and holomorphic functions $F_{k}$, predefined
on formal sectors $Q_{k}$, such that $F_{k+1} = \alpha\circ F_{k}$ at the
common boundary of adjacent sectors, $\ 0\leq k\leq n-1$ (cf. Theorem A 
below). 
As we note in the present 
article, this crucial symmetry property is also held (up to local gauge 
correction) by any contact form for which $\alpha$ is an irrational rotation
(as in the case of the irrational ellipsoids), so that a $\psi$ of charge $n$
is again determined by the holomorphic data $F_{k}$:

\vspace{.1in}

\noindent{\bf Theorem A.} {\em Suppose, for the ensemble} $(M,\lambda, j)$ {\em
that} ${\mathcal L}_{X}j = 0$ {\em and the return map} $\alpha$ {\em induced 
locally by the Reeb flow near 
a locally recurrent orbit} ${\mathcal P}$ 
{\em corresponds to an irrational rotation. For positive integer} $n$ 
{\em let} $Q_{k}$ {\em denote the formal sectors of a disc} 
$D\subset{\Bbb C}$ {\em defined by} 
\[2\pi\frac{k}{n} < \arg(z) < 2\pi\frac{k+1}{n}, \ \hspace{.2in} 0\leq k\leq 
n-1,\]
{\em with holomorphic functions} $w = F_{k}(z)$ {\em defined on} $Q_{k}$ 
{\em such that}
\[F_{k}\mid_{\overline{Q}_{k}\cap\overline{Q}_{k+1}} = \alpha
\circ F_{k+1}\mid_{\overline{Q}_{k}\cap\overline{Q}_{k+1}},\;\;\; \ \lim_{z
\rightarrow
0}F_{k}(z) = 0  \]   
{\em for each} $k$. {\em Modulo a local gauge correction of the form}
\[\hat{\lambda} = \lambda - df,\]
{\em for some smooth function} $f$, {\em defined in a tubular neighbourhood 
of} ${\mathcal P}$  {\em such that} $X_{\lambda}(f) = 0$, {\em these 
holomorphic data determine a finite--energy} $J$--{\em holomorphic curve of 
charge} $n$, {\em asymptotic to} ${\mathcal P}$.

\vspace{.1in}
  
But to how large a class of contact structures on ${\Bbb S}^{3}$ do these
correspondences apply? More broadly, to how large a class of structures do
the dynamic convexity results of Hofer, Wysocki and Zehnder apply? An important class
of examples of triples $(M,\la,j)$ is given by the contact form and almost 
complex endomorphism $j$ naturally induced on the three-manifold corresponding
to the unit tangent bundle of a Riemann surface $\Sigma$ by the geodesic flow
of a metric $g$ on $\Sigma$ (cf., e.g., \cite{Pat}). If the additional Lie 
symmetry
condition ${\mathcal L}_{X}j = 0$ is assumed then $g$ must have
Gaussian curvature $K\equiv 1$ and up to isometry we just have $\Sigma = {\Bbb S}^{2}$. But
the class of required structures is much broader if, instead of a Riemannian 
metric, one considers a {\em Finsler structure} on the two-sphere, i.e., a
hypersurface $M\subset T{\Bbb S}^{2}$ and a surjective submersion $\pi:M
\rightarrow{\Bbb S}^{2}$ such that for all $p\in{\Bbb S}^{2}$ the fibre
$\pi^{-1}(p)$ is a smooth, closed strictly convex curve enclosing the origin in
$T_{p}{\Bbb S}^{2}$ (cf. section two). The associated structure equations
derived from a canonical framing of $T^{*}M$ depend on three functions $I,J,
K$ over $M$. When the Finsler structure corresponds specifically to a 
Riemannian metric, it follows that $I\equiv 0$ and $K$ corresponds to the 
pullback of Gaussian curvature from $\Sigma$. The induced almost complex 
endomorphism $j$ on a Finsler contact manifold $(M,\lambda)$ satisfies 
${\mathcal L}_{X}j = 0$ precisely when $K\equiv 1$ for any $I,J$ (cf. Lemma \ref{fesy}). 

A famous class of contact structures due to A. Katok \cite{K} arises from
convex hypersurfaces $M\subset T{\Bbb S}^{2}$ on which $K\equiv 1$ and
the return map $\alpha$ associated with the Reeb flow near either of its
two distinguished periodic orbits is a rotation. The Finsler structures connected with
these examples incorporate what are known as {\em Randers metrics}, coming
from special perturbations of the norm derived from a Riemannian metric. 
When a Finsler structure on $\Sigma$ with $K\equiv 1$ admits a non-trivial 
Killing field, one may perturb to a family of non-Randers examples with the
same properties, as is shown in section 3.2, following independent 
observations of P. Foulon and W. Ziller \cite{Z}. All the hypersurfaces in
question are doubly covered by ${\Bbb S}^{3}$, such that the contact forms 
lift to
tight contact structures on the three-sphere. In sections three and four
a precise description of this lifting reveals that the Katok examples
correspond exactly with those induced by the restriction of $\lambda_{0}$
to ${\Bbb E}_{p,q}\subset{\Bbb C}^{2}$. It is then worth noting that the
non-Randers structures of section 3.2 lie genuinely beyond the standard class
of examples coming from convex hypersurfaces of ${\Bbb R}^{4}$. Returning
to the broader notion of dynamic convexity for contact structures on ${\Bbb
S}^{3}$, it is asked in the final section of this article whether a 
corresponding criterion can be found in terms of Finsler structures on
${\Bbb S}^{2}$. Theorem B provides this criterion for a given Finsler
structure, via a lower bound on the length $\ell$ of its shortest geodesic loop
(we shall say that the Finsler metric is dynamically convex if its associated 
contact form is dynamically convex):

\vspace{.1in}

\noindent{\bf Theorem B.} {\em Let} $F$ {\em be a Finsler metric on} ${\Bbb S}^2$ 
{\em such that} $K\geq\delta>0$.
{\em If} $\ell>\pi/\sqrt{\delta}$, {\em then} $F$ {\em is dynamically convex.}

\vspace{.1in}

In particular, by the results in \cite{HWZ} any such Finsler metric has 
either two or infinitely many closed geodesics.  We note that recently,
V. Bangert and Y. Long \cite{BL} have shown that {\it any} Finsler metric 
on ${\Bbb S}^2$ has two closed geodesics.

\bigskip

{\it Acknowledgements:} The authors would like to express their warm thanks to Z. Shen for
his helpful communications during their research towards this article.
We are also very grateful to V. Bangert and H.B. Rademacher for their comments and remarks.

\vspace{.1in}
                 
\section{$J$-holomorphic cylinders near locally recurrent orbits of a 
contact three-manifold}

\vspace{.1in}

Let $M$ denote a compact, oriented three-manifold with contact form $\lambda$
and associated plane-field $\xi\subset TM$ corresponding to ${\rm ker}(\lambda)$.
Let $X_{\lambda}$ denote the Reeb vector field associated with this structure
on $M$, together with an almost complex structure $J$ acting on $\xi$ such that
the symmetric tensor defined by $d\lambda(*,J*)\mid_{\xi}$ is positive 
definite. 
Consider a periodic orbit of the Reeb flow, denoted ${\mathcal P}$,
and a tubular neighbourhood $T_{\mathcal P}\subset M$.
If $\Delta$ represents a disc centred at the origin 
in ${\Bbb R}^{2}$, let $\tilde{\Delta}\subset M$ be an embedded image
such that the origin is mapped to the unique element $p_{0}$
of ${\mathcal P}\cap\tilde{\Delta}$, with $\tilde{\Delta}$ itself 
corresponding to a 
transverse slice of $T_{\mathcal P}$. The Reeb flow will be assumed moreover
to be Lyapunov--stable near ${\mathcal P}$ in the sense that for all 
$p\in\tilde{\Delta'}$,
where $\Delta'\subseteq\Delta$ is a sufficiently small disc centred at
the origin, there exists a unique solution 
$\gamma_{p}:[0,\infty)\to M$ to the equation
\[\frac{d\gamma_{p}}{dt} = X_{\lambda}(\gamma_{p}(t)) \ , \ \hspace{.2in}
\gamma_{p}(0) = p \ ,\]
which depends smoothly on both $t$ and $p$, and remains inside $T_{\mathcal P}$
for all $t\geq 0$. Given $p\in\tilde{\Delta'}$, we 
will define $(i) 
\ \tau(p)$ to be the smallest $t>0$ such that $\gamma_{p}(t)\in
\tilde{\Delta} \ , \ (ii) \ \Gamma_{p} := \gamma_{p}((0,\tau(p)])$
 \ and for each connected open neighbourhood of the origin, $\Omega
\subseteq\tilde{\Delta}$ , 
\[(iii) \ \Gamma(\Omega) := \cup_{p\in\Omega}
\Gamma_{p} \ .\]
We may now consider a recursively defined system of neighbourhoods
$\{\Omega_{k}\}$, such that $\Omega_{0}:=\tilde{\Delta'}$, while $\Omega_{k}$
denotes the origin--component of $\Gamma(\Omega_{k-1})\cap\Omega_{k-1}$.
The set $\Omega_{\infty}:=\cap_{k=0}^{\infty}\Omega_{k}$ was seen in \cite{HW}
proposition 1 to be conformally equivalent to a disc whenever it corresponds
to an open subset of $\tilde{\Delta}'$.
As in \cite{HW},  the Reeb flow will be said to be 
 ``locally recurrent'' near a periodic orbit ${\mathcal P}$ if it is 
Lyapunov--stable within a tubular neighbourhood $T_{\mathcal P}$
and for any sufficiently small 
embedded disc $\tilde{\Delta}$, corresponding to a transversal slice 
through $T_{\mathcal P}$ at some point $p_{0}$, the 
limit set $\Omega_{\infty}\subseteq\tilde{\Delta'}\subseteq\tilde
{\Delta}$ is open. An
orbit ${\mathcal P}$ itself may also be referred to as ``locally 
recurrent'' in this context. In passing we note that if the orbit is elliptic
of twist type, then by Moser's twist map theorem it will be
locally recurrent. Moreover, the local recurrence implies that 
the characteristic multipliers are unimodular. 

Under the assumption that the 
Reeb flow is locally recurrent near ${\mathcal P}$, we now select 
$\Omega_{\infty}\times\{\vartheta_{0}\}$
as coordinate disc within the initial Martinet tube (in which $\lambda
 = f\cdot(d\vartheta +xdy)$, where $\vartheta$ denotes the periodic 
coordinate and $f$ denotes a function such that $f(0,0,\vartheta)
\equiv\tau_{0}$ and $\nabla f(0,0,\vartheta)\equiv{\bold 0}$, though these 
facts are not used here). Without loss of 
generality, let $\vartheta_{0}$ be zero and 
consider the cylinder $\Omega_{\infty}\times[0,2\pi]$, which maps to the 
tube via the obvious identification ${\rm mod}(2\pi)$. The cylinder has 
${\mathcal P}$ 
as its axis, $x=y=0$, and the Reeb vector field in Martinet coordinates 
already looks like $\frac{1}{\tau_{0}}\frac{\partial}{\partial\vartheta}$ 
when restricted to ${\mathcal P}$. There is
no consequent loss of generality if we ``normalise'' $\lambda$ by the
constant multiple $\frac{1}{\tau_{0}}$, so that the minimal period is
effectively $1$. In particular this will ensure that the notion of 
``charge'' for a $J$-holomorphic mapping $\psi$ of the punctured disc, as 
introduced below, is consistent with the standard definition corresponding
to the limit of the integral of $\psi^{*}\lambda$ over the circle $|z| = c$
as $c$ approaches zero (and is specifically an integer). By analogy with
the standard construction of Darboux coordinates, the next step is to 
define 
\[{\mathcal C}:=\{(p,t) \ | \ p=(x,y)\in\Omega_{\infty} \ , \ \hspace{
.1in} 0\leq t\leq\tau(p) \ \} \ ,\]
and a homeomorphism 
\[h:{\mathcal C}\to\Omega_{\infty}\times[0,2\pi] \ , \hspace{.1in}
h\mid_{\Omega_{\infty}\times\{0\}} = 1 \ ,\]
which is smooth for all $0<t<\tau(p)$, coming from solutions of the
ordinary differential equation
\[\frac{d\gamma_{p}}{dt} = X_{\lambda}(\gamma_{p}(t)) \ .\]
It follows that on the interior of ${\mathcal C}$, the standard contact form
$\lambda_{0}$ and $\lambda':= h^{*}\lambda$ have the same Reeb vector field, 
corresponding to $\frac{\partial}{\partial t}$. 
We now consider the Cauchy--Riemann system
\[\pi((h^{-1}\psi)_{\eta}) + J\pi((h^{-1}\psi)_{\zeta}) = 0 \ 
\hspace{.1in}(*) \ ,\]
\[\lambda'((h^{-1}\psi)_{\zeta}) = -a_{\eta}, \ \hspace{.1in} \ 
\lambda'((h^{-1}\psi)_{\eta}) = a_{\zeta} \ \hspace{.1in}
(\dag) \ ,\]
satisfied by some finite--energy $J$--holomorphic map $\psi$ of a punctured 
neighbourhood
$D\setminus\{0\}\subset{\Bbb C}$ into the Martinet tubular neighbourhood of 
${\mathcal P}$ ($\times{\Bbb R}$). 
For a sufficiently ``thin'' neighbourhood of
${\mathcal P}$, the standard projection $(v_{1},v_{2},v_{3})\mapsto (v_{1},
v_{2})$ determines a linear isomorphism $\mu$ between $\xi' := 
ker(\lambda')$ and ${\Bbb R}^{2}$. Hence we define a $2\times 2$
matrix--valued function $j(x,y) = \mu\circ J\circ\mu^{-1}$, such that
${\bold x}:=(x,y)$ implies (*) can be written in the form
\[{\bold x}_{\eta}(z) + j{\bold x}_{\zeta}(z) = 0 \ .\]
Let $\alpha$ denote the diffeomorphism of $\Omega_{\infty}
\times\{0\}$ defined by the return map $\alpha(p):= 
\gamma_{p}(\tau(p))$, hence $\alpha(0) = 0$.
As seen in \cite{HW},  if ${\mathcal L}_{X_{\lambda}}J = 0$, then in a neighbourhood
of $0\in\Delta$, the smooth automorphism $\alpha$ is equivalent to a
rotation, via a
diffeomorphism $\varphi:\Delta''\rightarrow U\subseteq\Omega_{\infty}$ such
that $\varphi_{*}^{-1}\circ j\circ \varphi_{*} = j_{0}$ (the standard multiplication
by ${\bf i}$) . More specifically, let $\Omega_{\infty}'$
denote the simply connected domain inside $U$ which is stabilised by the Reeb flow.
The diffeomorphism $\hat{\alpha}:= \varphi^{-1}\circ\alpha
\circ\varphi$ \  then acts on $\varphi^{-1}(\Omega_{\infty}')\subseteq
\Delta''$ as an automorphism such that $\hat{\alpha}(0) = 0$
under the assumption of local recurrence. The additional assumption 
${\mathcal L}_{X_{\lambda}}J = 0$ implies that $\alpha^{*}j = j$, hence in
particular $\hat{\alpha}j_{0} = j_{0}\hat{\alpha}$, i.e., $\hat{\alpha}$
is a conformal automorphism. Modulo a conformal transformation
identifying $\varphi^{-1}(\Omega_{\infty}')$ with a disc, $\hat{\alpha}$
is then equivalent to a rotation. 
 
Now $\tilde{\bold x}_{\zeta}:=\varphi_{*}^{-1}({\bold x}_{\zeta})$ and 
$\tilde{\bold x}_{\eta}:=\varphi_{*}^{-1}({\bold x}_{\eta})$ implies
\[\tilde{\bold x}_{\eta}(\eta,\zeta) + j_{0}\cdot
\tilde{\bold x}_{\zeta}(\eta,\zeta) = {\bold 0} \ (\dag *) .\]
As described in \cite{HW}, each ``branch'' of
\[\Psi:=(h\circ(\varphi\times 1))^{-1}\psi\]
is defined smoothly in the interior and continuously up to the 
boundaries of a quasi--sector $Q_{k}$
in $D\setminus\{0\}$, with discontinuities arising at points $z_{0}$
lying on the smooth arcs, corresponding to $\vartheta^{-1}(0)$,that 
bound adjacent sectors (in the usual
way ``$\pm$'' will be used to denote opposite sides of the boundary). 
Discontinuities of the transverse 
projection of $\Psi$ are therefore described by the relations
\[\lim_{z\to z_{0}^{\pm}}\tilde{{\bold x}}(z):= \tilde{{\bold x}}^{\pm}
(z_{0}) \ \Rightarrow \ \hat{\alpha}(\tilde{{\bold x}}^{-}(z_{0})) = 
\tilde{{\bold x}}^{+}(z_{0}) \ .\]
Hence on each $Q_{k}\subset D\setminus\{0\} \ , \  (\dag *)$ defines a
holomorphic 
function  $w=F_{k}(z)$ which partially describes a branch of $\Psi$, 
such that
\[F_{k}\mid_{\overline{Q}_{k}\cap\overline{Q}_{k+1}}= \hat{\alpha}
\circ F_{k+1}\mid_
{\overline{Q}_{k}\cap\overline{Q}_{k+1}} \ ,\ \hspace{.1in} 0\leq 
k\leq n-1 \ ,\]
where $n$ denotes the charge (or asymptotic degree) of $\psi$ in 
relation to the periodic orbit ${\mathcal P}$.
 
Returning now to the particular form of the equations $(\dag)$ in
${\mathcal C}$ note that any $\lambda'$ with $\frac{\partial}{\partial t}$
as its Reeb vector field must take the general form 
\[\lambda' = dt + f_{1}({\bold x})dx + f_{2}({\bold x})dy \ .\]
Such a local presentation of the contact form allows decoupling of
$(\dag)$ into an inhomogeneous Cauchy--Riemann equation. This 
property is preserved under the diffeomorphism $\varphi$, however,
if it is assumed that ${\mathcal L}_{X_{\lambda}}
J = 0$, hence in particular the matrix $j$ above is independent of $t$.   
In this case, letting $\Re$ and $\Im$ signify real and imaginary parts,
and $u = t+{\bold i}a$, $(\dag)$ becomes

\[\frac{\partial}{\partial\bar{z}}[u\mid_{Q_{k}}] = -\{f_{1}(F_{k}(z))
\frac{\partial\Re(F_{k})}{\partial\bar{z}} + f_{2}(F_{k}(z))
\frac{\partial\Im(F_{k})}{\partial\bar{z}}\}\]
\[ = -\frac{1}{2}[(f_{1} +{\bold i}f_{2})
\circ F_{k}(z)]\cdot\overline{F'_{k}(z)} \ , \ \]
keeping in mind that this equation is defined smoothly
only on the interior of each quasi--sector $Q_{k}$. Define
$\omega := \frac{1}{2}(f_{1}+{\bold i}f_{2})(w,\bar{w})d\bar{w}$, so that
\[\lambda' = dt + 2\Re(\omega) \ , \ \hbox{and}\]
\[\frac{1}{2}[(f_{1} +{\bold i}f_{2})\circ F_{k}(z)]\cdot\overline{F'_{k}(z)}
d\bar{z} = F_{k}^{*}\omega \ .\]
Now ${\mathcal L}_{X_{\lambda}}\lambda = 0$ implies $\hat{\alpha}^{*}
\Re(\omega) = \Re(\omega)$. In particular, ${\bf f}:= (f_{1},f_{2})$
implies $\Re(\omega) = ({\bf f},*)$ with respect to the standard inner
product on ${\Bbb R}^{2}$, and hence $\hat{\alpha}_{*}^{t}{\bf f} = 
{\bf f}$. Similarly $\Im(\omega) = (j_{0}{\bf f},*)$, while
${\mathcal L}_{X_{\lambda}}J = 0$ implies $\hat{\alpha}_{*}^{t}j_{0} = 
j_{0}\hat{\alpha}_{*}^{t}$, so that
\[\hat{\alpha}^{*}\Im(\omega) = 
(\hat{\alpha}_{*}^{t}j_{0}{\bf f},*) = (j_{0}{\bf f},*) = \Im(\omega) \ .\]
It follows that $\hat{\alpha}^{*}\omega = \omega$, and hence
\[F_{k}^{*}\omega\mid_{\overline{Q}_{k}\cap\overline{Q}_{k+1}} = 
F_{k}^{*}(\hat{\alpha}^{*}\omega)\mid_{\overline{Q}_{k}\cap\overline
{Q}_{k+1}} = 
(\hat{\alpha}\circ F_{k})^{*}\omega\mid_{\overline{Q}_{k}\cap
\overline{Q}_{k+1}} \] 
\[ = F_{k+1}^{*}\omega\mid_{\overline{Q}_{k}\cap\overline{Q}_{k+1}} \ .\]
There now exists a 
continuous function $G(z,\bar{z})$ on $D$ such that
\[G(z,\bar{z})d\bar{z}\mid_{Q_{k}}:= F_{k}^{*}\omega \ , \ \hspace{.1in}
0\leq k\leq n-1 \ , \]
which is moreover continuously differentiable (cf. \cite{HW}), and hence
\[\hat{G}(z,\bar{z}) := \int_{D}\frac{G(\mu,\bar{\mu})}{\mu-z}
d\mu\wedge d\bar{\mu}\]
is a twice-continuously differentiable function on $D$. Moreover, the 
inhomogeneous equation satisfied by $u$ now has the form
\[\frac{\partial}{\partial\bar{z}}[u\mid_{Q_{k}}] = G, \;\;\; \ 0\leq k\leq n-1
 \ ,\] and with a little additional argument we have 

\begin{theorem}[\cite{HW}] Let $(\psi,a):D\setminus\{0\}\to M\times{\Bbb R}$ 
be a $J$--holomorphic curve of finite energy and charge $n$ at
 $z=0$,
asymptotic to a locally recurrent periodic orbit ${\mathcal P}$, near which 
${\mathcal L}_{X_{\lambda}}J = 0$. Consider any tubular neighbourhood
of ${\mathcal P}$ in $M$, diffeomorphic to $\Delta\times{\Bbb S}^{1}$ 
such that $\{0\}\times{\Bbb S}^{1}\approx{\mathcal P}$. There exists a diffeomorphic
change of coordinates in $\Delta\times[0,2\pi)$ such that on each 
quasi--sector $Q_{k}\subset D\setminus\{0\}$ the map $(\psi,a)$ 
can be expressed in the form 
\[(F_{k}(z), H_{k}(z) - \frac{1}{2\pi{\bold i}}\hat{G}(z,\bar{z}))\ , \ 
\hspace{.1in} 0\leq k\leq n-1 ,\]
where $F_{k} \ , \ H_{k}$ are holomorphic on $Q_{k}$ and 
continuous on $\overline{Q}_{k}$, such that
\[F_{k}\mid_{\overline{Q}_{k}\cap\overline{Q}_{k+1}} = \hat{\alpha}
\circ F_{k+1}\mid_{\overline{Q}_{k}\cap\overline{Q}_{k+1}} \ , \]
while each $H_{k}$ corresponds to an
analytic branch of $\frac{1}{2\pi{\bold i}}\log(\rho) \ , \ 
ord_{0}(\rho) = n$.
Moreover, the function $\hat{G}$ belongs to $C^{2}(D)$ and 
is bounded by $K|z|$. Finally, if $\alpha = 1$, then each
$F_{k}$ is the restriction of a single function $F$ holomorphic 
on $D$, \ $F(0) = 0$.
\end{theorem}

Conversely, given an ensemble $(M, \lambda, J)$, ${\mathcal L}_{X_{\lambda}}J 
= 0$,
with locally recurrent periodic orbit ${\mathcal P}$ corresponding to the 
asymptotic 
limit of some $J$--holomorphically embedded cylinder of finite energy, we may 
ask for the essential holomorphic data which determine such 
$J$-holomorphic curves in general near ${\mathcal P}$. Our main result in 
this direction is the following: 

\medskip

\noindent{\bf Theorem A.} {\it Suppose, for the ensemble $(M,\lambda, J)$ above, that the 
return
map $\alpha$ induced locally by the Reeb flow near a recurrent orbit 
${\mathcal P}$ 
corresponds to an irrational rotation. For positive integer $n$ let $Q_{k}$ 
denote
the formal sectors of a disc $D\subset{\Bbb C}$ defined by 
\[2\pi\frac{k}{n} < \arg(z) < 2\pi\frac{k+1}{n} \ \hspace{.1in} 0\leq k\leq 
n-1 \ ,\]
with holomorphic functions $w = F_{k}(z)$ defined on $Q_{k}$ such that
\[F_{k}\mid_{\overline{Q}_{k}\cap\overline{Q}_{k+1}} = \hat{\alpha}
\circ F_{k+1}\mid_{\overline{Q}_{k}\cap\overline{Q}_{k+1}} \ , \ \lim_{z
\rightarrow 0}F_{k}(z) = 0  \]   
for each $k$. Modulo a local gauge correction of the form
\[\hat{\lambda} = \lambda - df \ ,\]
for a smooth function $f$, defined in a tubular neighbourhood of 
${\mathcal P}$  
such that $X_{\lambda}(f) = 0$, these holomorphic data determine a 
finite--energy $J$--holomorphic curve of charge $n$, asymptotic to 
${\mathcal P}$.}

\medskip

\begin{proof} Returning to the general form of $\lambda$ with respect to 
the coordinate tube ${\mathcal C}$ above, we have
\[\lambda' = dt + f_{1}({\bf  x})dx + f_{2}({\bf x})dy = dt + \Re(\omega) \ ,\]
noting that $\hat{\alpha}^{*}\Re(\omega) = \Re(\omega)$, in the case of 
$\alpha$ an
irrational rotation (hence in particular the orbit of any point $w$ under 
$\alpha$ is 
dense in a circle of radius $|w|$). Now 
\[\Re(\omega) = \eta(\rho)d\nu + \zeta(\rho)d\rho = \eta(|w|)(xdy-ydx) + df 
\ ,\]
where $w = x + {\bf i}y$ \ , \ $\rho = |w|$ \ , \ $\nu = \arg(w)$ \ , \ and
\[f(\rho) = \int_{0}^{\rho}\zeta(r)dr \ .\]
Letting $\hat{\eta}(\rho^2) = \eta(\rho)$, and after making the gauge 
correction
$\hat{\lambda} := \lambda - df := dt + \Re(\hat{\omega})$, we have
\[\hat{\omega} = \hat{\eta}(w\bar{w})\cdot{\bf i}wd\bar{w} \ ,\]
and hence
\[\Phi(|w|^2) := \int_{0}^{|w|^2}\hat{\eta}(r)dr\]
implies $\hat{\omega} = {\bf i}\bar{\partial}\Phi(|w|^2)$. It should be 
noted at once
that $d\hat{\lambda} = d\lambda \ , \ X_{\lambda}(f) = 0$ together imply that 
$X_{\lambda}$
remains the Reeb vector field of the contact form defined by $\hat{\lambda}$. 
Recall that
\[F_{k}\mid_{\overline{Q}_{k}\cap\overline{Q}_{k+1}} = \hat{\alpha}
\circ F_{k+1}\mid_{\overline{Q}_{k}\cap\overline{Q}_{k+1}} \ , \ \lim_{z
\rightarrow
0}F_{k}(z) = 0  \]
indicates that there exists a single smooth function $|F|$ on $D$, 
corresponding to $F_{k}$
on each $Q_{k}$, hence we may define
\[G(z,\bar{z}) := {\bf i}\bar{\partial}\Phi(|F|^2) \ , \ \hbox{and} \ 
\hat{G}(z,\bar{z}) := 2\pi{\bf i}\cdot{\bf i}\Phi(|F|^2) = -2\pi\Phi(|F|^2) 
\ .\]   
Under the assumption that ${\mathcal P}$ is already the asymptotic limit of 
some finite--energy mapping of the punctured disc (or cylinder), we recall 
moreover from the discussion of \cite{HW}, that the return time $\tau(p)$ 
for each $p\in\Omega_{\infty}'$ is
constant (i.e., normalised to value 1). Note that the rotational symmetry of
the contact form ultimately implies that $\hat{G}$ is a real--valued 
function, compatible with the formal specification of simple sectors $Q_{k}$ 
as domains of the holomorphic functions $F_{k}$, and moreover that the 
$t$--component
of a $J$--holomorphic mapping derived from these data can be defined by $t = 
\frac{1}{2\pi}\arg(z^{n})$. Now from \cite{HW},  Theorem 1, we may write 
\[u = t +{\bf i}a = \frac{1}{2\pi{\bf i}}(\log(z^{n}) + 2\pi\Phi(|F|^2)) \ ,\]
from which it follows that the associated $a$--component must be
\[a(z) = \frac{-1}{2\pi}\log(|z|^{n}) - \Phi(|F|^2) \ .\]
One or two remarks should be made concerning the almost complex structure 
with respect
to which the mapping $\Psi$, determined sector-wise by the $F_{k}(z)$ and 
the corresponding
analytic branches of $\frac{1}{2\pi{\bf i}}(\log(z^{n}) + 2\pi\Phi(|F|^2))$, 
may be said
to represent a $J$--holomorphic mapping of charge $n$ of the cylinder into 
${\mathcal T}_{{\mathcal P}}\times{\Bbb R}$, 
for a tubular neighbourhood ${\mathcal T}_{{\mathcal P}}$. It is 
straightforward to see
that in a sufficiently thin tube, the standard projection mapping defines an 
isomorphism
between ${\Bbb R}^2$ and the contact planes of both $\lambda$ and $\hat
{\lambda}$ at any
given point which we will denote by $\beta:\xi\rightarrow\hat{\xi}$ , i.e., 
$\beta({\bf v}) = {\bf v} - \hat{\lambda}({\bf v})X$ \ , where
$X$ denotes the common Reeb vector field of $\lambda$ and $\hat{\lambda}$.
The almost complex structure $\hat{J} := \beta\circ J\circ
\beta^{-1}$ is then automatically induced on ${\rm ker}(\hat{\lambda})$. Moreover, 
given 
${\bf v'}\in {\rm ker}(\hat{\lambda})$, we have
\[d\hat{\lambda}({\bf v'},\hat{J}\cdot{\bf v'}) = d\hat{\lambda}(\beta({\bf 
v}),\hat{J}\cdot
\beta({\bf v})) = d\hat{\lambda}(\beta({\bf v}),\beta(J\cdot{\bf v}))\]
\[ = d\hat{\lambda}({\bf v} - \hat{\lambda}({\bf v})X,J\cdot{\bf v} - \hat{
\lambda}
(J\cdot{\bf v})X) = d\lambda({\bf v},J\cdot{\bf v}) \ .\]
Hence the quadratic form $d\hat{\lambda}(*,\hat{J}*)\mid_{\hat{\xi}}$ is also 
positive definite. 

It remains now to check that the pseudoholomorphic curve $\Psi$, defined with 
respect to 
$\hat{\lambda} \ , \ \hat{J}$ in ${\mathcal T}_{\mathcal P}\times{\Bbb R}$ is 
of finite energy.
Following \cite{HWZ}, let ${\mathcal F}$ denote the space of smooth functions 
$h:{\Bbb R}\rightarrow
[0,1]$ such that $h'\geq 0$, and define extensions $\hat{\lambda}_{h}$ of the 
contact
form from ${\mathcal T}_{{\mathcal P}}$ to ${\mathcal T}_{{\mathcal P}}\times
{\Bbb R}$ such that
$\hat{\lambda}_{h}(p,a) := h(a)\cdot\hat{\lambda}(p)$. The ``energy'' of 
$\Psi$ is
then defined as 
\[E(\Psi) := \sup_{{\mathcal F}}\int_{D\setminus\{0\}}\Psi^{*}d\hat{\lambda}_
{h} \ .\]
Clearly,
\[\int_{D\setminus\{0\}}\Psi^{*}d\hat{\lambda}_{h} = \int_{|z| = 1}\Psi^{*}
\hat{\lambda}_{h} - 
\lim_{\varepsilon\rightarrow 0}
\int_{|z| = \varepsilon}\Psi^{*}\hat{\lambda}_{h}\]
\[ =  \int_{|z| = 1}h(a)\Psi^{*}\hat{\lambda} - \lim_{\varepsilon\rightarrow 0}
\int_{|z| = \varepsilon}h(a)\Psi^{*}\hat{\lambda} \ .\]
Now 
\[\Psi^{*}\hat{\lambda} = -\eta(|F|)d\nu(F) + \frac{n}{2\pi}d\arg(z) \ ,\]
where the formula
\[\nu = \arg(w) = \tan^{-1}\left(-{\bf i}\frac{w-\bar{w}}{w+\bar{w}}\right)\]
implies 
\[d\nu(F) = 2\Re(\frac{d\nu}{dw}\cdot F'(z)dz) = \Re(\frac{-{\bf i}F'}{F}dz) 
\ .\]
(Here we have adopted a harmless abuse of notation, in the sense that the 
ratios
$\frac{F_{k}'}{F_{k}}$ define a single continuous function on $D\setminus
\{0\}$ when $\alpha$ is a rotation).  
Note, inside any $|z| = \varepsilon$ sufficiently small, that $z = 0$ is the 
unique zero of the smooth
function $|F(z)|$, and hence
\[\left| \int_{|z| = \varepsilon}\Psi^{*}\hat{\lambda} - n\right| \leq 
\int_{|z| = \varepsilon}\frac{|\eta(|F|)||F'|}
{|F|}\cdot\varepsilon\cdot d\arg(z)  \ , \]
where the latter quantity approaches zero as $\varepsilon$ goes to zero,
if it is recalled that
$\eta(\rho)$ is a smooth function such that $\eta(0) = 0$ . Now $|z|
\rightarrow 0$ implies $a\rightarrow +\infty$, and thus
\[E(\Psi) = \sup_{{\mathcal F}}\left(\int_{|z| = 1}h(-\Phi(|F(z)|^2))
\cdot\Psi^{*}\hat{\lambda} - n\cdot h(+\infty)\right)
 \ ,\]
which is clearly finite.

\end{proof}
 
In the next section we will consider a class of examples of contact 
structures to which the
above theorem may be applied, but it should first be estalished that 
holomorphic data of the
sort specified are in plentiful supply. For arbitrary $n$, consider
$f(z)$ holomorphic on $D\subset{\Bbb C}$, such that $ord_{0}(f) > n$,
and an arbitrary irrational rotation of the form $\alpha = e^{2\pi{\bf i}c} 
\ , \ 0 < c < 1$.
Let $D$ be divided into formal sectors as in theorem 1, with 
\[F_{0}(z) := z^{-nc}f(z) \ , \ F_{k+1}(z) = \alpha\cdot F_{k}(z)  \ , 
\ 0\leq k\leq n-1,\]
noting that $F_{n}(z) = F_{0}(z)$. For convenience we may define 
$z^{-nc} = e^{-ncLog(z)}$, where $Log(z)$ denotes an analytic branch 
of the complex logarithm defined on ${\Bbb C}$ minus the positive real
axis.

\vspace{.1in}

\section{Examples}

\subsection{Canonical coframing}

Let $\Sigma$ be a closed oriented connected surface.
A smooth {\it Finsler structure} on $\Sigma$ is
 a smooth hypersurface $M\subset T\Sigma$
for which the canonical projection $\pi:M\to\Sigma$ is a surjective
submersion having the property that for each $x\in \Sigma$, the $\pi$-fibre
 $\pi^{-1}(x)=M\cap T_{x}\Sigma$ is a smooth, closed, strictly
 convex curve enclosing
the origin $0_x\in T_{x}\Sigma$.

Given such a structure it is possible to define a canonical coframing
$(\a,\b,\c)$ on $M$ that satisfies the following structural
 equations (see \cite[Chapter 4]{BCZ}):
\begin{align}
d\a&=-\b\wedge\c,\label{eq1}\\
d\b&=-\c\wedge(\a-I\b),\label{eq2}\\
d\c&=-(K\a-J\c)\wedge\b.\label{eq3}
\end{align}
where $I$, $K$ and $J$ are smooth functions on $M$.
The function $I$ is called the {\it main scalar}
of the structure and it vanishes if and only if
$M$ is the unit circle bundle of a Riemannian metric.
When $I=0$, i.e. when the Finsler structure is Riemannian,
$K$ is the $\pi$-pullback of the Gaussian curvature.

Let $X_1$, $X_2$ and $X_3$ be the vector fields on $M$ that are
dual to the coframing $(\a,\b,\c)$. 
The form $\a$ is the canonical contact form of $M$ whose Reeb vector
field is the geodesic vector field $X_1$.

As a consequence of (\ref{eq1}--\ref{eq3}) the framing $(X_1,X_2,X_3)$
satisfies the commutation relations:
\begin{equation}\label{comm}
[X_3,X_1]=X_2,\quad [X_2,X_3]=X_1+IX_2+JX_3,\quad [X_1,X_2]=KX_3.
\end{equation}

Note that the coframing also defines a natural almost complex structure
$\mathbb J$ on $\xi=\mbox{\rm ker}\a$. Indeed, we may set:
\[{\mathbb J}(xX_2+yX_3)=yX_2-xX_3.\]
If we let $\eta=xX_2+yX_3\in \mbox{\rm ker}\a$ then using (\ref{eq1}) we see that:
\[d\a(\eta,{\mathbb J}\eta)=x^2+y^2\]
and thus ${\mathbb J}$ is compatible with the contact structure.

\begin{lemma} ${\mathcal L}_{X_1}{\mathbb J}=0$ if and only if $K=1$.
\label{fesy}
\end{lemma}

\begin{proof} Let $\phi_t$ be the flow of $X_1$. Note that ${\mathcal L}_
{X_1}{\mathbb J}=0$ if and only if $d\phi_t$ is an isometry of the inner 
product in 
$\mbox{\rm ker}\a$, given $d\a(*,{\mathbb J}*)$.
(Recall that $\phi_t$ preserves $d\a$.)

Let $\eta=xX_2+yX_3\in \mbox{\rm ker}\a$ and write:
\[d\phi_{t}(\eta)=x(t)X_2+y(t)X_3.\]
Thus the flow $\phi_t$ is an isometry of the inner product $d\a(*,{\mathbb J}
*)$
if and only if
\[\frac{d}{dt}(x(t)^2+y(t)^2)=0,\]
equivalently if and only if
\begin{equation}
x\dot{x}+y\dot{y}=0.
\label{iso}
\end{equation}
Write
\[\eta=x(t)d\phi_{-t}(X_2)+y(t)d\phi_{-t}(X_3)\]
and differentiate with respect to $t$ to obtain
\[0=\dot{x}X_2+x[X_1,X_2]+\dot{y}X_3+y[X_1,X_3].\]
Using the structure equations (\ref{comm}) and regrouping we have:
\[0=(\dot{x}-y)X_2+(\dot{y}+xK)X_3,\]
hence
\[\dot{x}=y,\]
\[\dot{y}+Kx=0.\]
If we use the last two equations in (\ref{iso}) we see that
$\phi_t$ is an isometry if and only if
\[xy(1-K)=0\]
and thus ${\mathcal L}_{X_1}{\mathbb J}=0$ if and only if $K=1$.
\end{proof}

\begin{remark}{\rm Note that the proof above shows that in general, if we
let $\eta\in\mbox{\rm ker}\a$ and write
\[d\phi_{t}(\eta)=x(t)X_2+y(t)X_3\]
then
\[\dot{x}=y,\]
\[\dot{y}+Kx=0.\]
We will use this fact later on.}
\end{remark}

We now recall some global consequences of $K=1$ as explained by R. Bryant 
in \cite{Br}.
The first thing to observe is that the structure equations imply:

\begin{align}
&\phi_t^*\a=\a,\\
\label{ec1}&\phi_{t}^*\b=\cos t\, \b+ \sin t\,\c\\
\label{ec2}&\phi_{t}^*\c=-\sin t \,\b+ \cos t\,\c
\end{align}

Suppose that $\Sigma$ is geodesically complete and connected. Then it can be 
shown that
$\Sigma$ is diffeomorphic to $S^2$ and there exists a unique orientation 
reversing isometry
$A:\Sigma\to\Sigma$ such that $dA|_{M}=\phi_{\pi}$ (we call $A$ a {\it 
quasi-antipodal map}). Moreover
for any point $p\in \Sigma$ every unit speed geodesic leaving $p$ passes 
through $A(p)$
at distance $\pi$ and $\Sigma$ has diameter $\pi$. According to \cite
[Proposition 4]{Br} we have the following dichotomy:

\begin{enumerate}
\item $A^2$ is the identity on $\Sigma$ in which case all geodesics are 
closed with the same
minimal period $2\pi$;
\item $A^2$ has exactly two fixed points, say $\sigma$ and $A(\sigma)$. 
Moreover there exists
a positive definite inner product on $T_{\sigma}\Sigma$ that is preserved 
by $d(A^2)(\sigma):T_{\sigma}\Sigma\to T_{\sigma}\Sigma$
and there is an angle $\theta \in (0,2\pi)$ such that $d(A^2)(\sigma)$ is 
counterclockwise rotation by $\theta$ in this inner product.

\end{enumerate}

Here we will be mostly interested in case (2). We note that this case has 
two possible subcases.
Suppose that $\theta/2\pi$ is rational and write $\theta=2\pi p/q$ where
$0<p\leq q$ with $p$ and $q$ coprime. Then $A^{2q}$ is the identity and thus
$\phi_{2\pi q}=\mbox{\rm identity}$, that is, every orbit of $\phi$ is closed 
with period $2\pi q$, although some orbits may have smaller minimal period.

When $\theta/2\pi$ is irrational, then the iterates of $A^2$ are dense in a 
circle of isometries of the Finsler surface $\Sigma$. This circle of 
isometries fixes $\sigma$ and $A(\sigma)$
and $\Sigma$ is rotationally symmetric about $\sigma$. Hence the geodesic flow
of $\Sigma$ is completely integrable with a ``Clairaut'' first integral.
The surface is also 
symmetric with respect
to $A$ about a circle (the equator) $E$. The unit tangent vectors to $E$ 
determine two closed
orbits $\gamma_{\pm}$. The equator divides $\Sigma$ into two disks 
$D_{\sigma}$ and $D_{A(\sigma)}$
which contain $\sigma$ and $A(\sigma)$ respectively. ($A$ maps $D_{\sigma}$ 
to $D_{A(\sigma)}$ and fixes $E$
setwise.)

Let $S$ be the subset of $M$ given by those pairs $(x,v)$ where $x\in E$ 
and $v$
points inside the region $D_\sigma$. The set $S$ (diffeomorphic to $E\times 
(0,\pi)$)
is a section of the geodesic flow with return map $\phi_{2\pi}|_{S}=d(A^2)
|_{S}$.
We see that there are no other closed orbits besides $\gamma_{\pm}$.

If $\pi\circ\gamma_{+}$ is the closed geodesic which travel around $E$ 
counterclockwise
(as seen from $\sigma$) then we find the lengths of $\pi\circ \gamma_{\pm}$ 
to be:
\[\ell_{+}:=\ell(\pi\circ\gamma_{+})=2\pi-d(p,A^2(p)),\]
\[\ell_{-}:=\ell(\pi\circ\gamma_{-})=2\pi+d(A^{2}(p),p),\]
where $p$ is any point in $E$. Since $\theta/2\pi$ is irrational we see that
$\ell_{\pm}/2\pi$ are irrational and in view of (\ref{ec1}) and (\ref{ec2}) 
we conclude that
$\gamma_{\pm}$ are elliptic orbits of $\phi$.
It is also clear that they are locally recurrent.

\subsection{Examples with only two closed geodesics}
\label{newex}

Summarizing the discussion above, Theorem A can be applied to a Finsler 
metric on ${\Bbb S}^2$ with $K=1$
and $\theta/2\pi$ irrational. Examples of such metrics are given by the well
known Katok examples \cite{K} analyzed by W. Ziller in \cite{Z}.
The fact that these metrics have $K=1$ is proved by Z. Shen \cite{S} (see
also \cite[Section 5]{R} for a discussion of these examples).

However these are not the only examples and a Katok type construction also 
gives a larger class as we now explain.

Suppose $F$ is a Finsler metric on ${\Bbb S}^2$ with $K=1$ and $A^2=Id$ 
(i.e. all the geodesics
are closed and with the same minimal period). Suppose in addition that 
$F$ admits a nontrivial
Killing field $V$. We may suppose without loss of generality that the flow 
$f_t$ of $V$ is such that
$f_{2\pi}=Id$.

 Define a 1-parameter family of Finsler metrics $G_{\varepsilon}$ 
($\varepsilon$
small) by giving its co-metric $G^*_{\varepsilon}$ in
$T^*{\Bbb S}^2$ as follows:
\begin{equation}
G^*_{\varepsilon}(x,p)=F^*(x,p)+\varepsilon p(V(x)).
\label{newfin}
\end{equation}
An unpublished result of P. Foulon asserts that $G_{\varepsilon}$ also 
has $K=1$ 
(this can be checked along the lines of the calculations in \cite{Z}).
The Katok examples arise when $F$ is the standard Riemannian metric with 
$K=1$ in ${\Bbb S}^2$.

We now note that Bryant \cite{Br1,Br2} has produced several families of 
Finsler metrics
with $K=1$ and $A^2=Id$. Among them there are subfamilies with rotational 
symmetry which are not of Randers type. For example Theorem 10 in \cite{Br1} 
gives a 1-parameter family
of rotationally invariant projectively flat Finsler metrics with $K=1$.

To complete our construction we note that Ziller in \cite{Z} observed that
for a Finsler metric defined by (\ref{newfin}) with $\varepsilon$ irrational, 
the only
closed geodesics of $G_{\varepsilon}$ are those invariant under $f_t$.
Hence there will be only two closed geodesics for $\varepsilon$ irrational.

We remark that the results in \cite{BRZ} show that the Katok examples are 
the only examples (up to isometry, of course)
of {\it Randers metrics} with $K=1$. A Randers metric is a Finsler metric of 
the
form $\sqrt{g_{x}(v,v)}+\theta_{x}(v)$, where $g$ is a Riemannian metric and 
$\theta$ is a 1-form.

\section{Lifting Finsler metrics to contact structures on ${\Bbb S}^3$}

Let $F$ be a Finsler metric on ${\Bbb S}^2$. The Lagrangian $\frac{1}{2} F^2$ 
gives rise to a Legendre
transform $\ell_{F}:T{\Bbb S}^{2}\setminus\{0\}\to T^{*}{\Bbb S}^{2}
\setminus\{0\}$ and 
if we let $\la$ be the Liouville 1-form
on $T^*{\Bbb S}^{2}$, it is well known that $\ell_{F}^{*}\la$ restricted 
to the unit sphere bundle
$M$ of $F$ is precisely the contact form $\omega_{1}$ from the previous 
section.

On $T^{*}{\Bbb S}^2$ we also have a Finsler co-metric $F^*$ such that 
$F=F^*\circ\ell_{F}$ and a corresponding
unit co-sphere bundle $M^*$.
We consider on ${\Bbb S}^2$ the canonical metric with curvature $1$. It 
has an associated Legendre tranform
$\ell_{0}:T{\Bbb S}^2\to T^{*}{\Bbb S}^2$. In what follows $|\cdot|$
denotes the norm of the canonical metric in both $T{\Bbb S}^2$ and $T^*
{\Bbb S}^2$.

Observe that there exists a unique smooth function $g:T^{*}{\Bbb S}^2
\setminus\{0\}\to {\mathbb R}^+$
such that $F^*(x,p)=g(x,p)|p|$. ($g$ is homogeneous of degree zero in $p$.)

Let $r:M^{*}_{0}\to M^*$ be given by 
\[r(x,p)=(x,p/g(x,p)).\]

We identify ${\Bbb S}^3$ with ${\Bbb S}{\Bbb U}(2)$ and ${\Bbb S}^2$ with 
the matrices in 
${\Bbb S}{\Bbb U}(2)$ of the form
\[x=\left(\begin{array}{cc}
it&z\\
-\bar{z}&-it\\
\end{array}\right),\]
where $t\in\mathbb R$, $z\in \mathbb C$ and $t^2+|z|^2=1$. If $A\in 
{\Bbb S}{\Bbb U}(2)={\Bbb S}^3$ and $x\in {\Bbb S}^2$ is as above
then $A^{-1}xA\in {\Bbb S}^2$. This also gives a natural embedding of 
${\Bbb S}^2$ in $\mathbb R^3$. 
Below we will often write $(t,z)$ to indicate a point in ${\Bbb S}^2$ 
instead of the corresponding matrix in ${\Bbb S}{\Bbb U}(2)$.

 Consider the matrices:
\[j=\left(\begin{array}{cc}
0&1\\
-1&0\\
\end{array}\right)\in {\Bbb S}^2,\;\;\;\;\;\
k=\left(\begin{array}{cc}
0&i\\
i&0\\
\end{array}\right)\in {\Bbb S}^2. \]

Now let ${\mathbb G}:{\Bbb S}^3\to M_{0}$ be the following map:

\[\mathbb G(A)=(A^{-1}j A,A^{-1}kA).\]
It is easy to see that the unit vectors $A^{-1}j A$, $A^{-1}kA$ in 
$\mathbb R^3$ are orthogonal.
Let $A\in {\Bbb S}{\Bbb U}(2)$ be written as:
\[A=\left(\begin{array}{cc}
w_1&w_2\\
-\bar{w}_{2}&\bar{w}_{1}\\
\end{array}\right)\]
where $(w_1,w_2)\in\mathbb C^2$ and $|w_1|^2+|w_2|^2=1$. The map $\mathbb G$ 
may also be written as
\begin{equation}
\mathbb G(A)=[(2\Im(w_{1}w_{2}),\bar{w}_{1}^2+w_{2}^2),\;(-2\Re(w_{1}w_{2}),
i(\bar{w}_{1}^2-w_{2}^2))].
\label{formg}
\end{equation}

If $\lambda_0$ is the canonical contact form on ${\Bbb S}^3$ ($\la_0(x)=
\langle ix,*\rangle$) and $\a^0$ is the contact form
of $M_0$, then it is not hard to check that (and it is explicitly done in 
\cite[Proposition 1.1]{CO}):
\begin{equation}
{\mathbb G}^*\a^0=2\lambda_0.
\label{gf}
\end{equation}
We now claim (compare with \cite[Proposition 1.1]{CO}):

\begin{lemma} Let $f:=1/g\circ\ell_{0}:T{\Bbb S}^2\setminus \{0\}\to 
\mathbb R^+$.
Then 
\[{\mathbb G}^* \ell_{0}^{*} r^*\left(\la|_{M^*}\right)=2(f\circ{\mathbb G})
\lambda_0\,.\]
\label{lift}
\end{lemma}

\begin{proof} Let $\tau:T^*{\Bbb S}^2\to{\Bbb S}^2$ be the canonical 
projection and note that
\[r^*\la_{(x,p)}(\xi)=p/g(d\tau(dr(\xi)))=\frac{1}{g}\la_{(x,p)}(d\tau(\xi))\]
and that
\[\ell_{0}^*(\la/g)=\omega_{0}^1/(g\circ\ell_{0}).\]
Combining these two equalities with (\ref{gf}) the lemma follows.

\end{proof}

The lemma is saying that the geodesic flow of a Finsler metric on ${\Bbb S}
^2$ is 
(up to a double covering) smoothly conjugate to the 
Reeb flow of a (tight) contact form on ${\Bbb S}^3$ of the form 
$h\,\la_{0}$ where 
$h=2(f\circ\mathbb G)$ and $f$ is related to
the Finsler metric as described in the lemma. Conversely, if we have a contact 
form $h\,\la_{0}$ with
$h(-A)=h(A)$, then $h$ will give rise to functions $f$ and $g$ as above. In 
general, for an arbitrary $h$ invariant under the 
antipodal map, the hypersurface $M^*$
of $T^*{\Bbb S}^2$ determined by $g$ does not need to be fibrewise strictly 
convex, but it will be clearly starshaped.
Note that in the proof of the lemma we did not really need $M^*$ to come 
from a Finsler metric.
If $M^*$ is just starshaped, then $\la|_{M^*}$ is also a contact 
form.

Summarizing, the lemma gives this: there is a 1-1 correspondence between
starshaped hypersurfaces of $T^*{\Bbb S}^2$ and smooth positive functions 
$h$ on ${\Bbb S}^3$ with $h(A)=h(-A)$.
If in addition $M^*$ is fibrewise strictly convex we obtain a Finsler metric.

\section{Ellipsoids and the Katok examples}

In contact geometry there is a well studied class of examples given by the 
ellipsoids
\[\mathbb E_{p,q}:=\{(w_1,w_2)\in \mathbb C^2:\;\;p|w_1|^2+q|w_{2}|^2=1\},\]
where $p$ and $q$ are positive real numbers. The restriction of $\la_0$ to
$\mathbb E_{p,q}$ determines a Reeb flow whose dynamics is very simple: the 
flow is just
$\phi_{t}(w_1,w_2)=(w_1 e^{ipt},w_{2} e^{iqt})$. There are two periodic 
orbits corresponding to
$w_1=0$ and $w_2=0$. These are the only periodic orbits if $p/q$ is 
irrational, whereas $\mathbb E_{p,q}$
is foliated by periodic orbits if $p/q$ is rational (but not all will have 
the same minimal period).

On the other hand, in Finsler geometry there are the well known Katok 
examples, which in particular provide examples
of Finsler metrics with only two closed geodesics. The Finsler co-metric 
of the Katok examples in geodesic polar
coordinates $(r,\phi)\in (0,\pi)\times [0,2\pi]$ is
\begin{equation}
F_{\varepsilon}^*(r,\phi,p_{r},p_{\phi})=\sqrt{p_{r}^2+\frac{1}{\sin ^2 r}
p_{\phi}^2}+\varepsilon\, p_{\phi},
\label{ka}
\end{equation}
where $\varepsilon\in (-1,1)$. We will show in this section that the 
ellipsoids and the Katok examples
are related precisely by the correspondence described in the previous section.

Let $h:{\Bbb S}^3\to\mathbb R^+$ be the function
\begin{equation}
h(w_1,w_2)=\frac{1}{p|w_1|^2+q|w_2|^2}.
\label{forh}
\end{equation}
Also let $\varphi:{\Bbb S}^3\to \mathbb E_{p,q}$ be given by
$\varphi(w_1,w_2)=\sqrt{h(w_1,w_2)}\,(w_1,w_2)$. It is easy to check that
\[\varphi^*\la_{0}=h\,\la_{0}\]
and hence we can think of the Reeb flow of the ellipsoids as being defined 
on ${\Bbb S}^3$ with contact form
$h\,\la_0$ where $h$ is given by (\ref{forh}).
Clearly $h(A)=h(-A)$. To find $f$ we write $\mathbb G(A)=(x,v)\in M_{0}$
where $x=(t,z)$ and $v=(b,\eta)$. The expressions for $t,z,b,\eta$ are given 
by (\ref{formg}). Using them we derive
\[\bar{w}_{1}^2=\frac{z-i\eta}{2},\]
\[w_{2}^2=\frac{z+i\eta}{2}.\]
Using that $(w_1,w_2)\in {\Bbb S}^3$ we obtain
\[f(x,v)=\frac{1}{(p-q)|z-i\eta|+2q}.\]
We now introduce geodesic polar coordinates $(r,\phi)$ on ${\Bbb S}^2$ such 
that
\begin{align*}
&t=\cos r,\\
&z=\sin r e^{i\phi}.
\end{align*}
Thus $(b,\eta)=(\dot{t},\dot{z})$. In these coordinates the Legendre transform
$\ell_{0}$ is simply:
\[p_{r}=\dot{r},\;\;\;\;p_{\phi}=\dot{\phi}\sin^2 r.\]
We can now compute $g=1/(f\circ\ell_{0}^{-1})$ in $(r,\phi,p_{r},p_{\phi})$
-coordinates:
\begin{align*}
g(r,\phi,p_{r},p_{\phi})&=2q+(p-q)\sqrt{\dot{r}^2\cos^2 r+\sin^2 r\,(1+\dot{
\phi})^2}\\
&=2q+(p-q)\sqrt{p_{r}^2\cos^2 r+\sin^2 r+2p_{\phi}+\frac{p_{\phi}^2}{\sin^2 
r}}.
\end{align*} 
This is the expression of $g$ on $M^*_{0}$, that is, when in addition
\[p_{r}^2+\frac{p_{\phi}^2}{\sin^2 r}=1.\]
Thus if we simplify it further we obtain
\[g(r,\phi,p_{r},p_{\phi})=2q+(p-q)(1+p_{\phi})=p+q+(p-q)\,p_{\phi}.\]
But in view of (\ref{ka}) the value of $F^*_{\varepsilon}$ on $M^*_{0}$
is just $1+\varepsilon p_{\phi}$. Note that by homogeneity a Finsler metric is 
completely determined by its value on $M^*_{0}$.
Hence if we choose $p$ and $q$ such that $p+q=1$, the ellipsoids induce, 
under the correspondence described in the previous section, exactly the Katok 
examples with $\varepsilon=p-q$.

\begin{remark} {\rm In \cite[Theorem 1.1]{HWZ} the authors show that a 
strictly convex hypersurface
$S\subset \mathbb R^4$ carries either 2 or infinitely many periodic orbits. 
The ellipsoids, of course, provide examples
of this dichotomy and the authors point out that it is not true that the 
first alternative holds
only for the irrational ellipsoids. They also remark in \cite[Page 200]{HWZ} 
that M. Herman constructed examples
of hypersurfaces $S$ which are $C^\infty$-close to the ellipsoid, admit only 
two periodic orbits but have a transitive flow on $S$. Apparently these 
examples of Herman are unpublished, but we remark here that it is
easy to construct such examples using the results above combined with Katok's 
main result in \cite{K}.
Indeed, Katok shows in \cite{K} that given any $r$, one can approximate 
$F_{\varepsilon}$ in the $C^r$ topology 
by a Finsler metric $F$ with ergodic geodesic flow and only two closed 
geodesics.
We have shown above that under the lifting procedure described in Lemma 
\ref{lift}, $F_{\varepsilon}$ lifts
to an ellipsoid, and hence $F$ will give rise to a smooth hypersurface $S$ 
which is $C^r$-close to an ellipsoid
and such that the flow on $S$ has only two closed orbits and is transitive.

Note that the examples described in Subsection \ref{newex} give rise to a 
new 1-parameter
family of hypersurfaces in $\mathbb R^4$ exhibiting the same dynamics and 
dichotomy as the ellipsoids.

Finally we note that if under the lifting described in Lemma \ref{lift} a 
Finsler metric gives rise to
a strictly convex hypersurface, then it will carry either 2 or infinitely 
many closed geodesics. Recently,
V. Bangert and Y. Long \cite{BL} have shown that {\it any} Finsler metric 
on ${\Bbb S}^2$ has two closed geodesics.
In the next section we will describe geometric conditions on the Finsler 
metric that ensure 
{\it dynamical convexity} so the main results in \cite{HWZ} can be applied.}

\end{remark}

\section{Dynamically convex Finsler metrics}

We first recall the definition of {\it dynamically convex} contact form $\la$ on a closed, connected
and oriented manifold $M$ with $\pi_2(M)=0$.

The Reeb vector field $X$ is transversal to the contact structure $\xi=\mbox{\rm ker}\,\la$ so that naturally
we have a splitting $T_{x}M=\mathbb R X(x)\oplus \xi_{x}$. Let $\phi_t$ be the flow of $X$.
Clearly $d\phi_{t}(x):\xi_{x}\to \xi_{\phi_{t}(x)}$ is symplectic with respect to $d\la$.
A contractible periodic solution $x$ with period $T$ has an integer-valued 
index $\mu(x,T)$ which we briefly recall
(cf. \cite{HWZ}).

Let $D\subset \mathbb C$ be the closed unit disk. Choose a smooth map $\sigma:
D\to M$
such that $\sigma(e^{2\pi {\bf i} t/T})=x(t)$. Choose a {\it symplectic} 
trivialization $\psi:\sigma^*\xi\to D\times\mathbb R^2$
of the symplectic bundle $\sigma^*\xi$ with symplectic form $\sigma^* d\la$.
Here $\mathbb R^2$ is endowed with the standard symplectic form.
We can use this trivialization to define a symplectic arc
\[\Phi:[0,T]\to Sp(1)\]
by setting
\[\Phi(t)=\psi(e^{2\pi {\bf i} t/T})\circ d\phi_{t}(x(0))\circ\psi^{-1}(1).\] 
Any such symplectic arc has a {\it Conley-Zehnder index} $\mu(\Phi)$, whose 
definition we recall below,
and we define
\[\mu(x,T)=\mu(\Phi).\]
It can be checked that $\mu(x,T)$ is well defined and does not depend on the 
choices made.

There are several possible ways to define $\mu(\Phi)$. We shall present
the definition that is most appropriate to our purposes. Our reference for 
what follows is \cite[Section 3]{HWZ}.

Let $\Phi:[0,T]\to Sp(1)$ be a smooth arc with $\Phi(0)=I$ and set
$A(t):=-J\dot{\Phi}(t)\Phi(t)^{-1}$. Then $A(t)$ is a smooth path of symmetric 
matrices and $\Phi$ solves the linear differential equation:
\[\dot{\Phi}=JA\Phi,\;\;\;\;\Phi(0)=I\]
for $t\in[0,T]$.
Given $\tau\in\mathbb R$, suppose $v$ (not identically zero) solves the first 
order differential equation
\begin{equation}
-J\dot{v}(t)-A(t)v(t)=\tau\,v(t)
\label{eigen}
\end{equation}
and satisfies the periodic boundary condition $v(0)=v(T)$.
Then $v$ is an eigenvector with eigenvalue $\tau$ of the operator 
\[L_{A}(v)=-J\dot{v}-A(t)v\]
defined on the space of $T$-periodic $H^1$-maps into $\mathbb R^2$.
Since $v$ never vanishes we may
choose a smooth angle
$\varphi(t)$ such that
\[e^{2\pi {\bf i} \varphi(t)}=\frac{v(t)}{|v(t)|},\;\;\;t\in [0,T].\]
Define the winding number $\Delta(\tau,A)\in\mathbb Z$ by
\[\Delta(\tau,A):=\varphi(T)-\varphi(0).\]
One can see that $\Delta(\tau,A)$ depends only on the eigenvalue $\tau$
and not on the eigenfunction $v$. For every integer $k\in\mathbb Z$
there are precisely two eigenvalues (counting multiplicities) $\tau_1$ and 
$\tau_2$ such that $k=\Delta(\tau_1,A)=\Delta(\tau_2,A)$.
Hence we may label the eigenvalues $\tau$ of $L_A$ by their winding numbers. 
Indeed, set $\Delta(\tau_{k},A)=[\frac{k}{2}]$ for
$k\in\mathbb Z$ and $\tau_{k}\leq\tau_j$ for $k\leq j$.
With this labeling we define the
{\it Conley-Zehnder index} of $\Phi$ as
\[\mu(\Phi):=\max\{k:\;\tau_{k}<0\}.\]

%\[\mu(\Phi):=\left\{\begin{array}{ll}
%2k+1,&\;\;\;\;\mbox{\rm if}\;I(\Phi)\subset (k,k+1)\\
%2k,&\;\;\;\;\mbox{\rm if}\;k\in I(\Phi)\\
%\end{array}\right. \]
%where $k\in\mathbb Z$.

A contact form $\la$ on $M$ is said to be {\it dynamically convex} if 
$\mu(x,T)\geq 3$ for every contractible
periodic solution of the Reeb vector field $X$.
In \cite{HWZ} it is shown that if a contact form $\la$ on $\mathbb S^3$ comes from
a strictly convex hypersurface in $\mathbb R^4$, then it is dynamically convex.

Recall from the previous sections that a Finsler metric on ${\Bbb S}^2$ 
defines a contact form in ${\Bbb S}{\Bbb O}(3)$.
We shall say that the Finsler metric is dynamically convex if its associated 
contact form is dynamically convex.

Given a Finsler metric on ${\Bbb S}^2$, let $\ell$ be the length of the 
shortest geodesic loop.

\medskip

\noindent{\bf Theorem B.} {\it Let $F$ be a Finsler metric on ${\Bbb S}^2$ such that 
$K\geq\delta>0$.
If $\ell>\pi/\sqrt{\delta}$, then $F$ is dynamically convex.
In particular by the results in \cite{HWZ} any such Finsler metric has either 
two or infinitely many closed geodesics.}

\medskip

\begin{proof} By a simple rescaling argument we may assume that
$K\geq 1$ and that $\ell>\pi$.

Consider a closed geodesic $\gamma$ with length $T$.
The vector fields $X_2$ and $X_3$ provide a trivialization of $\xi$ without 
the need of extending 
the bundle to $D$. Given $\eta\in \xi$, write $\eta=aX_2+bX_3$. Then
\[\eta\mapsto (b,a)\]
gives a symplectic trivialization between $(\xi,d\omega_1)$ and $\mathbb R^2$ 
with the canonical symplectic form.
(Note that we need to swap $a$ and $b$ since $d\omega_1=-\omega_2\wedge 
\omega_3$.)
Recall that if we write $d\phi_{t}(\eta)=x X_2+y X_3$, then $\dot{x}=y$ and 
$\ddot{x}+Kx=0$.
Thus, using the trivialization, we see that $\Phi$ satisfies the linear 
differential equation:
\begin{equation*}
\dot{\Phi}=JA\Phi
\label{linear}
\end{equation*}
where
\[J=\left(\begin{array}{cc}

0&-1\\
1&0\\

\end{array}\right)\;\;\;
\mbox{\rm and}\;\;\;
A=\left(\begin{array}{cc}

1&0\\
0&K\\

\end{array}\right). \]
Take $v$ an eigenvector of $L_A$ with eigenvalue $\tau$
and set $w(t):=\frac{v(t)}{|v(t)|}$.
Using (\ref{eigen}) it is straightforward to check that $w(t)$ satisfies
\[\dot{w}=J(A+\tau I)w-w\langle w,J(A+\tau I)w\rangle.\]
But $2\pi \dot{\varphi} Jw=\dot{w}$ and thus $2\pi\dot{\varphi}=\langle (A+
\tau I)w,w\rangle$.
Since $K\geq 1$, we see that $\langle (A+\tau I)w,w\rangle\geq 1+\tau$ for all 
$t\in [0,T]$ and hence
\[\Delta(\tau,A)=\varphi(T)-\varphi(0)\geq \frac{T}{2\pi}(1+\tau).\]
Since $\Delta(\tau_3,A)=1$ we deduce
\begin{equation}
\frac{T}{2\pi}(1+\tau_3)\leq 1.
\label{taui}
\end{equation}
Let us combine this inequality with the lower bound on $\ell$.
If $\gamma$ is a simple closed curve in
${\Bbb S}^2$, the orbit $t\mapsto x(t)=(\gamma(t),\dot{\gamma}(t))$ will not 
be contractible in $M$.
If $\gamma$ is not a simple closed geodesic, then it must have length $T\geq 
2\ell>2\pi$.
Inequality (\ref{taui}) shows that $\tau_3<0$
and by the definition of $\mu$
we see that $\mu(x,T)\geq 3$ as desired.

\end{proof}

\begin{remark}{\rm The proof of the theorem indicates that one cannot expect 
to get dynamical convexity
just assuming positive curvature. In fact if $\gamma$ is a closed geodesic
with length $T$ such that $K=1$ along it, then the proof of the theorem shows 
that
\begin{equation}
1=\frac{T}{2\pi}(1+\tau_3).
\label{k=1}
\end{equation}
Now consider a convex surface of revolution whose equator
has $K=1$. It is clear that one can make the length $T$ of the equator as 
short as desired at the expense
of increasing the curvature elsewhere. But using (\ref{k=1}) 
we see that the 
Conley-Zehnder index of the equator iterated twice (to get a contractible 
curve in the unit sphere
bundle) is $\leq 2$ if $T$ is small enough since $\tau_3\geq 0$ in this case.
A concrete example of this is given by the ellipsoid in $\mathbb R^3$:
\[\frac{x^2+y^2}{a^2}+z^2=1.\]
The equator $z=0$ has curvature $K=1$ and length $2\pi a$, so for $a=1/2$ we
get that the equator iterated twice has $\tau_3=0$ and $\mu=1$.
Note that in this case $1\leq K\leq 16$.
This example shows that Theorem B is in fact sharp.

Hence there are positively curved Riemannian metrics on ${\Bbb S}^2$ such 
that under the lifting procedure
described in Lemma \ref{lift} they do not give rise to convex hypersurfaces 
in $\mathbb R^4$.

The results of Hofer, Wysocki and Zhender in \cite{HWZ} treat a lot more than 
just the dichotomy ``two or infinitely many'' closed orbits. Their results
state that the Reeb flow has a closed disk as a surface of section whose 
boundary
is an unknotted closed orbit $P$ with $\mu=3$. Moreover, any other closed 
orbit has to be linked
with $P$.} 

\end{remark}

If the Finsler metric is {\it reversible}, then one can obtain a lower bound
for $\ell$ from an upper bound on curvature.
Indeed if $0<K<a$, then $\ell>2\pi/\sqrt{a}$. For Riemannian metrics
this is a classical result (cf. \cite[Theorem 3.4.8]{Kl}) and for reversible
Finsler metrics the proof is quite the same (cf. \cite[Theorem 4]{R}).

Thus Theorem B implies:

\begin{corollary} A strictly 1/4-pinched positively curved reversible Finsler
metric is dynamically convex.
\end{corollary}

The corollary seems new even for the case of Riemannian metrics.
Given a Finsler metric $F$, the {\it reversibility} $r$ of $F$ is
(cf. \cite{R}):
\[r:=\max_{(x,v)\in M}\{F(x,-v)\}\geq 1.\]
Clearly $r=1$ iff $F$ is reversible.
In \cite[Theorem 4]{R} Rademacher gives a lower bound for the length $L$ of 
the shortest closed geodesic as follows. If $0<K\leq 1$, then
\[L\geq \pi(1+1/r).\]

The same bound holds for $\ell$ \cite{Ra1} and hence
we may also conclude that a Finsler metric with
$\left(1-\frac{1}{1+r}\right)^2\leq K<1$ is dynamically convex.

\end{document}